\documentclass[12pt,reqno]{amsart}
\usepackage{pb-diagram} 
\usepackage{amsmath,amssymb}

\usepackage{amsthm}
\usepackage[colorlinks = true,
            linkcolor = blue,
            urlcolor  = blue,
            citecolor = blue,
            anchorcolor = blue]{hyperref}
\usepackage[active]{srcltx}
\usepackage{pb-diagram}

\numberwithin{equation}{section}

\newtheorem{theorem}{Theorem}[section]

{Proposition}[section]
\newtheorem{corollary}
{Corollary}[section]
\newtheorem{remark}{Remark}[section]

\newfont{\bb}{msbm10 at 12pt}

\def\<{\langle}     
\def\>{\rangle}

\def\Paps{{\bf P}_+}

\newcommand{\bal}{\begin{align}}      \newcommand{\eal}{\end{align}}
\newcommand{\ba}{\begin{array}}      \newcommand{\ea}{\end{array}}
\newcommand{\bc}{\begin{center}}     \newcommand{\ec}{\end{center}}
\newcommand{\be}{\begin{enumerate}}  \newcommand{\ee}{\end{enumerate}}
\newcommand{\beq}{\begin{eqnarray}}  \newcommand{\eeq}{\end{eqnarray}}
\newcommand{\beQ}{\begin{eqnarray*}} \newcommand{\eeQ}{\end{eqnarray*}}
\newcommand{\bi}{\begin{itemize}}    \newcommand{\ei}{\end{itemize}}
\newcommand{\bt}{\begin{tabular}}    \newcommand{\et}{\end{tabular}}
\newcommand{\bdm}{\begin{displaymath}} \newcommand{\edm}{\end{displaymath}}

\newcommand{\D}{D\!\!\!\!/\,}

\newcommand{\Ss}{R\!\!\!\!/\,}

\begin{document}

\title[Total mean curvature and first Dirac eigenvalue]{Total mean curvature and first Dirac eigenvalue}       
\author{Simon Raulot}
\address{Laboratoire de Math\'ematiques R. Salem
UMR $6085$ CNRS-Universit\'e de Rouen
Avenue de l'Universit\'e, BP.$12$
Technop\^ole du Madrillet
$76801$ Saint-\'Etienne-du-Rouvray, France.}
\email{\href{simon.raulot@univ-rouen.fr}{simon.raulot@univ-rouen.fr}}

\begin{abstract}
In this note, we prove an optimal upper bound for the first Dirac eigenvalue of some hypersurfaces in Euclidean space by combining a positive mass theorem and the construction of quasi-spherical metrics. As a direct consequence of this estimate, we obtain an asymptotic expansion for the first eigenvalue of the Dirac operator on large spheres in three dimensional asymptotically flat manifolds. We also study this expansion for small geodesic spheres in a three dimensional Riemannian manifold. We finally discuss how this method can be adapted to yield similar results in the hyperbolic space.
\end{abstract}



\date{\today}   

\maketitle 
\pagenumbering{arabic}

 
\section{Introduction}


A smooth and connected $n$-dimensional Riemannian manifold $(M,g)$ is asymptotically flat if there exists a compact set $K\subset M$ and a diffeomorphism $\Phi:M\setminus K\rightarrow\mathbb{R}^n\setminus B$, for a closed ball $B$, such that in the asymptotically flat coordinates $x=(x_1,\cdots,x_n)$ given by $\Phi$, we have
\begin{eqnarray}\label{MetricAF}
g_{ij}=\delta_{ij}+\sigma_{ij},
\end{eqnarray}
where $\sigma_{ij}$ are smooth functions satisfying
\begin{eqnarray*}
\sigma_{ij}=O(|x|^{-\tau}),\quad\partial_k\sigma_{ij}=O(|x|^{-\tau-1}),\quad\partial_k\partial_l\sigma_{ij}=O(|x|^{-\tau-2})
\end{eqnarray*}
as $|x|\rightarrow+\infty$ for some constant $\tau>\frac{n-2}{2}$ and for all $i,j,k,l=1,\cdots,n$. Moreover, we require that the scalar curvature $R$ of $g$ is integrable. In this situation, the ADM mass of $(M,g)$ is defined as 
\begin{eqnarray*}
m_{ADM}(M,g):=\frac{1}{2(n-1)\omega_{n-1}}\lim_{r\rightarrow+\infty}\sum_{i,j=1}^n\int_{S_r}\Big(\frac{\partial g_{ij}}{\partial x_j}-\frac{\partial g_{jj}}{\partial x_i}\Big)\nu_idS_r
\end{eqnarray*}
where $\omega_{n-1}$ is the area of the unit sphere $\mathbb{S}^{n-1}\subset\mathbb{R}^n$ and $S_r\subset M$ is a large coordinate sphere of radius $r$ with outward normal $\nu$. It is a well-known fact that the limit in the right-hand side above is finite and that its values does not depend on the chosen asymptotically flat coordinates so that it defines a geometric invariant of $(M,g)$ (see \cite{Bartnik1,Chrusciel1}). When $(M,g)$ is a space-like slice of a time-symmetric isolated gravitational system, the ADM mass measures the energy of this system. In this situation, the scalar curvature expresses the energy density of the system so that it is natural to conjecture that if $R\geq 0$ then $m_{ADM}(M,g)\geq 0$ with equality if, and only if, $(M,g)$ is isometric to the Euclidean space. This is the famous positive mass theorem which is known to be true by the works of Schoen and Yau \cite{SchoenYau5,SchoenYau3} and Witten \cite{Witten1}. 

If $n=3$, the positive mass theorem is equivalent to a result of Shi and Tam \cite{ShiTam1} which states that if $(\Omega_0,g)$ is a compact Riemannian manifold with nonnegative scalar curvature and whose boundary is a $2$-sphere with positive Gauss curvature and positive mean curvature $H$ then 
\begin{eqnarray}\label{ShiTam-Inequality}
\int_\Sigma Hd\Sigma\leq \int_\Sigma \mathcal{H}_0d\Sigma.
\end{eqnarray}
Here $\mathcal{H}_0$ denotes the mean curvature of the unique strictly convex isometric embedding of $(\Sigma,\gamma)$, $\gamma:=g_{|\Sigma}$, in the Euclidean space $(\mathbb{R}^{3},\delta)$ whose existence and unicity are ensured by the Weyl's embedding theorem \cite{Nirenberg,Pogorelov}. In a certain sense, the Shi-Tam's inequality localizes the positive mass theorem. The proof of the inequality (\ref{ShiTam-Inequality}) is based on a low regularity positive mass theorem which can be obtained using the Witten's method. This approach requires the existence of a spin structure on $M$ which allows to define spinors as well as the Dirac operator associated to the Riemannian metric $g$. One of the key point in Witten's proof is that both the scalar curvature and the mass appear in the corresponding Bochner formula, the first as the curvature term in its pointwise version and the second as the boundary-at-infinity contribution in its integral version. 

In this note, we remark that the Shi-Tam's approach can be used to derive a new upper bound for the first nonnegative eigenvalue of the extrinsic Dirac operator of some compact hypersurfaces in Euclidean space. More precisely, we will prove:
\begin{theorem}\label{NewUpperBound}
Let $\Sigma$ be a compact, orientable, mean-convex and star-shaped hypersurface with positive scalar curvature in the $n$-dimensional Euclidean space. Then the first nonnegative eigenvalue $\lambda_1(\Sigma,\gamma)$ of the extrinsic Dirac operator $\D_\gamma$ of $(\Sigma,\gamma)$ satisfies
\begin{eqnarray}\label{NewUpperBound1}
\lambda_1(\Sigma,\gamma)\leq\frac{\int_\Sigma \mathcal{H}_0\,d\Sigma}{2|\Sigma|}.
\end{eqnarray}
Here $\mathcal{H}_0$ denotes the mean curvature function of $\Sigma$ in $\mathbb{R}^{n}$ and $|\Sigma|$ its volume. Equality occurs if, and only if, the hypersurface is a round sphere.
\end{theorem}
Note that $\lambda_1(\Sigma,\gamma)$ corresponds to the absolute value of the first eigenvalue of the intrinsic Dirac operator of $(\Sigma,\gamma)$ which has to be positive by the Friedrich inequality \cite{Friedrich3} since we assumed that its scalar curvature is positive. We refer to the monograph \cite{BourguignonHijaziMilhoratMoroianu} where the reader can find most of the proofs of the results we used in this paper concerning spin geometry. When $(\Sigma,\gamma)$ is a $2$-sphere with positive Gauss curvature, the Weyl's embedding theorem implies the following intrinsic upper bound:
\begin{corollary}\label{2SphereDirac}
Let $(\Sigma,\gamma)$ be a $2$-sphere with positive Gauss curvature then
\begin{eqnarray*}
\lambda_1(\Sigma,\gamma)\leq\frac{\int_\Sigma \mathcal{H}_0\,d\Sigma}{2|\Sigma|}.
\end{eqnarray*}
Equality occurs if, and only if, $(\Sigma,\gamma)$ is a round sphere.
\end{corollary}
The proof of Theorem \ref{NewUpperBound} is a direct combination of the construction of Bartnik \cite{Bartnik1} and Shi-Tam \cite{ShiTam1} of quasi-spherical metrics (and more generally its generalization by Eichmair, Miao and Wang \cite{EichmairMiaoWang}) and a positive mass theorem for manifolds with boundary due to Herzlich \cite{Herzlich1,Herzlich2}.

This new estimate leads to several applications. First, it allows to study the asymptotic behavior of the first nonnegative Dirac eigenvalue of large spheres in three dimensional asymptotically flat manifolds. More precisely, we have: 
\begin{theorem}\label{AsymptoticDiracSphere}
Let $(M,g)$ be an asymptotically flat $3$-dimensional manifold and let $S_r$ be a coordinate sphere of radius $r>0$ in some chart at infinity. Then if $\lambda_1(S_r,g_r)$ denotes the first nonnegative eigenvalue of the extrinsic Dirac operator on $S_r$ endowed with the metric $g_r$ induced by $g$, it holds that
\begin{eqnarray}\label{IntegralInequality1}
\lambda_1(S_r,g_r)|S_r|=\frac{1}{2}\int_{S_r} H_r\,dS_{r}+4\pi m_{ADM}(M,g)+o(1).
\end{eqnarray}
\end{theorem}
Here $H_r$ represents the mean curvature of $S_r$ in $(M,g)$ and $|S_r|$ the area of $(S_r,g_r)$. An immediate consequence of Theorem \ref{AsymptoticDiracSphere} and the positive mass theorem is:
\begin{corollary}\label{LS-Mass}
Let $(M,g)$ be an asymptotically flat $3$-manifold with nonnegative scalar curvature. Then
\begin{eqnarray*}
\lim_{r\rightarrow+\infty}\Big(\lambda_1(S_r,g_r)|S_r|-\frac{1}{2}\int_{S_r} H_r\,dS_{r}\Big)\geq 0,
\end{eqnarray*}
and equality holds if, and only if, $(M,g)$ is isometric to the Euclidean space $\mathbb{R}^3$. 
\end{corollary}
 
Another direct consequence of Corollary \ref{2SphereDirac} is provided by the small-sphere limit of the first eigenvalue of the Dirac operator. More precisely, we prove:
\begin{theorem}\label{SmallSphereExpansion}
Let $(M,g)$ be a $3$-dimensional Riemannian manifold, $p$ be a fixed interior point of $M$ and $S_r$ be the geodesic sphere of radius $r$ centered at $p$ with induced metric $\gamma_r$. For $r$ small enough, we have
\begin{eqnarray}\label{SSE}
\lambda_1(S_r,\gamma_r)=\lambda_1(S_r,\delta_r)+\frac{R(p)}{36}r+O(r^3)
\end{eqnarray}
where $\lambda_1(S_r,\delta_r)=1/r$ is the first eigenvalue of the Dirac operator on a $2$-dimensional round sphere of radius $r$.
\end{theorem}
This means that, at least for $n=3$, the scalar curvature at a point can be recovered from the first eigenvalues of the Dirac operator on small geodesic spheres. We also deduce immediately a result similar to Theorem \ref{AsymptoticDiracSphere} in this situation, namely: 
\begin{corollary}\label{SmallSphereDiracV}
Let $(M,g)$ be a Riemannian manifold of dimension $3$, $p$ be a fixed interior point of $M$ and $S_r$ be the geodesic sphere of radius $r$ centered at $p$ with induced metric $\gamma_r$. For $r$ small enough, we have
\begin{eqnarray}\label{SSE}
\lambda_1(S_r,\gamma_r)|S_r|=\frac{1}{2}\int_{S_r} H_r\,dS_{r}+\frac{\pi}{3}R(p)r^3+O(r^5).
\end{eqnarray}
In particular, if $R(p)\geq 0$, then 
\begin{eqnarray*}
\lim_{r\rightarrow 0}\frac{1}{r^3}\Big(\lambda_1(S_r,\gamma_r)|S_r|-\frac{1}{2}\int_{S_r} H_r\,dS_{r}\Big)\geq 0,
\end{eqnarray*}
and equality occurs if, and only if, $(M,g)$ is flat at $p$.
\end{corollary}

It turns out that this approach can easily be extended to the hyperbolic setting. More precisely, if $(\Sigma,\gamma)$ is a compact hypersurface isometrically embedded in the $n$-dimensional hyperbolic space $\mathbb{H}^n(-\kappa^2)$ with constant sectional curvature $-\kappa^2$ for some $\kappa>0$, we are naturally lead to consider the first nonnegative eigenvalues $\lambda_1^\pm(\Sigma,\gamma)$ of zero order modifications $\D_\gamma^\pm$ of the extrinsic Dirac operator $\D_\gamma$ (see Section \ref{HyperbolicPart} for precise definitions). Then we prove the following upper bound regarding these eigenvalues:
\begin{theorem}\label{NewUpperBoundHyperbolique}
Let $(\Sigma,\gamma)$ be a hypersurface isometrically embedded into $\mathbb{H}^n(-\kappa^2)$ with mean curvature $\mathcal{H}_0$ and which is homeomorphic to a $(n-1)$-sphere with sectional curvature $K>-\kappa^2$. Then the first nonnegative eigenvalue $\lambda_1^\pm(\Sigma,\gamma)$ of the Dirac-type operator $\D_\gamma^\pm$ satisfies
\begin{eqnarray}\label{NewUpperBoundHyp1}
\lambda_1^\pm(\Sigma,\gamma)\int_\Sigma \cosh(\kappa r)\,d\Sigma\leq\frac{1}{2}\int_\Sigma \mathcal{H}_0\cosh(\kappa r)\,d\Sigma.
\end{eqnarray}
Equality occurs if, and only if, the hypersurface is a geodesic sphere centered at the origin.
\end{theorem}
Here, without loss of generality, we assumed that $\Sigma$ encloses a region $\Omega$ which contains $o=(0,\cdots,0,1/\kappa)\in\mathbb{H}^n(-\kappa^2)\subset\mathbb{R}^{n,1}$ where $\mathbb{R}^{n,1}$ denotes the $(n+1)$-dimensional Minkowski space and $r$ is the geodesic distance of a point from $o$. Combining this estimate with a lower bound on $\lambda_1^\pm(\Sigma,\gamma)$ gives the following Minkowski-type inequality when $n=3$:
\begin{corollary}\label{MinkowskiHyperbolic}
Let $(\Sigma,\gamma)$ be a $2$-sphere with sectional curvature bounded from below by $-\kappa^2$ embedded into $\mathbb{H}^3(-\kappa^2)$ with mean curvature $\mathcal{H}_0$. Then 
\begin{eqnarray}\label{MinkowskiHyperbolic1}
\int_\Sigma\mathcal{H}_0\cos(\kappa r)\,d\Sigma\geq 4\sqrt{\frac{\pi}{|\Sigma|}+\frac{\kappa^2}{4}}\int_\Sigma\cosh(\kappa r)\,d\Sigma.
\end{eqnarray}
Moreover, equality occurs if, and only if, $(\Sigma,\gamma)$ is a geodesic sphere centered at the origin.
\end{corollary}
This last inequality should be compared to similar inequalities obtained by Ge, Wang and Wu \cite{GeWangWu} and by Brendle, Hung and Wang \cite{BrendleHungWang}. 

The paper is organized as follows. In Section \ref{UpperBoundDiracEuclidean}, we give the proof of Theorem \ref{NewUpperBound}. The first Dirac eigenvalue expansions on large spheres and small spheres are respectively proved in Section \ref{LargeSphereLimit} and Section \ref{SmallSphereLimit}. The last part is devoted to the study of the hyperbolic setting. 

 
\section{A new upper bound for the first Dirac eigenvalue}\label{UpperBoundDiracEuclidean}


In order to prove Theorem \ref{NewUpperBound}, we need to quickly review the proof of the Shi-Tam inequality (\ref{ShiTam-Inequality}). The first main ingredient in their proof is the construction of quasi-spherical metrics initiated by Bartnik \cite{Bartnik4}. Let $(\Sigma,\gamma)$ be a compact and strictly convex hypersurface embedded in the $n$-dimensional Euclidean space. Then, the set $\mathbb{R}^n\setminus\Omega$, where $\Omega$ is the compact domain enclosed by $\Sigma$ in $\mathbb{R}^n$, is foliated by $\Sigma_{\rho}$ the hypersurface at distance $\rho$ from $\Sigma$. In the following, we will identify $\mathbb{R}^n\setminus\Omega$ with $\mathcal{E}:=\Sigma\times[0,+\infty[$. Throughout this identification, the Euclidean metric on $\mathcal{E}$ can be expressed as $d\rho^2+\gamma_\rho$ where $\gamma_\rho$ is the induced metric on $\Sigma_\rho$ and $\gamma_0=\gamma$. Then it is proved in \cite{ShiTam1} that for any positive function $u_0$ on $\Sigma$, there exists a positive function $u$ on $\mathcal{E}$ such that the Riemannian metric 
\begin{eqnarray}\label{QS-metric}
g_u:=u^2d\rho^2+\gamma_\rho
\end{eqnarray}
is scalar flat with $u(\cdot,0)=u_0$. One can say even more since in fact $(\mathcal{E},g_u)$ is an asymptotically flat manifold with well-defined ADM mass. The second crucial fact in the proof of (\ref{ShiTam-Inequality}) is that the function
\begin{eqnarray*}
\rho\in[0,\infty)\mapsto\int_{\Sigma_\rho}\mathcal{H}_\rho(1-u^{-1})d\Sigma_\rho
\end{eqnarray*}
where $\mathcal{H}_\rho$ denotes the mean curvature of $\Sigma_\rho$ in $\mathbb{R}^n$, is nonincreasing in $\rho$ and tends to a multiple of the ADM mass of $(\mathcal{E},g_u)$ as $\rho$ goes to infinity. In particular, it holds that
\begin{eqnarray}\label{Monotone}
\int_{\Sigma}\mathcal{H}_0(1-u^{-1}_0)d\Sigma\geq c(n)m_{ADM}(\mathcal{E},g_u)
\end{eqnarray}
for some positive constant $c(n)$ depending only on $n$. So, as soon as one can show that the mass of $(\mathcal{E},g_u)$ is nonnegative, we get
\begin{eqnarray}\label{GeneralInequality}
\int_\Sigma \mathcal{H}_0d\Sigma\geq \int_\Sigma u^{-1}_0\mathcal{H}_0d\Sigma.
\end{eqnarray}
This can be done when $(\Sigma,\gamma)$ is a $2$-sphere with positive Gauss curvature bounding a compact Riemannian $3$-manifold $(\Omega_0,g)$ with nonnegative scalar curvature and positive mean curvature $H$. Indeed, since $(\mathcal{E},g_u)$ has a compact inner boundary which is isometric to $(\Sigma,\gamma)$ with mean curvature equals to $u_0^{-1}\mathcal{H}_0$, if we choose the initial value $u_0=\mathcal{H}_0/H$, one can glue $(\Omega_0,g)$ and $(\mathcal{E},g_u)$ along their common boundaries to get an asymptotically flat manifold for which the positive mass theorem holds. The inequality (\ref{GeneralInequality}) exactly yields the Shi-Tam inequality. 

In our situation, we will directly apply a positive mass theorem for asymptotically flat manifolds with compact inner boundary to $(\mathcal{E},g_u)$. This result, due to Herzlich \cite{Herzlich1,Herzlich2}, ensures that the ADM mass is nonnegative when the first eigenvalue of the Dirac operator of $(\Sigma,\gamma)$ satisfies a certain lower bound. More precisely, it states:
\begin{theorem}$($\cite{Herzlich1,Herzlich2}$)$\label{HerzlichPMT-n}
Let $(M,g)$ be a $n$-dimensional Riemannian spin asymptotically flat manifold with nonnegative scalar curvature and with a compact inner boundary $\Sigma$. Assume that the first nonnegative eigenvalue $\lambda_1(\Sigma,\gamma)$ is positive and satisfies
\begin{eqnarray}\label{DiracMeanCur-n}
\lambda_1(\Sigma,\gamma)\geq\frac{1}{2}H
\end{eqnarray}
where $H$ is the mean curvature of $\Sigma$ in $M$. Then the mass is nonnegative and if the mass is zero, $(M,g)$ is flat and the mean curvature is constant equals to $H=2\lambda_1(\Sigma,\gamma)$. 
\end{theorem}
This result is sharp since the exterior of round balls in Euclidean space are flat manifolds with zero mass for which (\ref{DiracMeanCur-n}) is an equality. In our conventions, the mean curvature of an $(n-1)$-dimensional round sphere with radius $r>0$ in Euclidean space is $(n-1)/r$. We are now in position to prove Theorem \ref{NewUpperBound}. 

\vspace{0.2cm}

{\it Proof of Theorem \ref{NewUpperBound}:} First assumed that $(\Sigma,\gamma)$ is a strictly convex hypersurface in $\mathbb{R}^n$. As recalled above, one can solve the quasi-spherical metric problem on $\mathcal{E}$ with the initial value 
\begin{eqnarray*}\label{InitialValueBV}
u_0=\frac{1}{2\lambda_1(\Sigma,\gamma)}\mathcal{H}_0>0.
\end{eqnarray*}
Then the metric $g_u$ defined by (\ref{QS-metric}) yields, on $\mathcal{E}$, an asymptotically flat metric with zero scalar curvature for which the inner compact boundary is isometric to $(\Sigma,\gamma)$ and has constant mean curvature equals to $H=2\lambda_1(\Sigma,\gamma)$. It is immediate to see that Theorem \ref{HerzlichPMT-n} applies so that we deduce from (\ref{GeneralInequality}) the inequality (\ref{NewUpperBound1}). If now we assume that equality occurs, then the ADM mass of $(\mathcal{E},g_u)$ is zero and the positive mass theorem implies that the metric $g_u$ is flat. From the Gauss formula, one can compute that the Riemann curvature tensor $Riem_u$ of $g_u$ is given by
\begin{eqnarray*}\label{RiemFlat}
Riem_u(e_i,e_j,e_i,e_j)=(1-u^{-2})Riem_\rho(e_i,e_j,e_i,e_j)
\end{eqnarray*}
where $Riem_\rho$ is the Riemann curvature tensor of $\Sigma_\rho$ and $\{e_i\,/\,1\leq i\leq n-1\}$ is a local orthonormal frame on $(\Sigma,\gamma)$ parallel translated in the direction of $\frac{\partial}{\partial \rho}$. Since the metric $g_u$ is flat and $\Sigma_\rho$ is strictly convex in $\mathbb{R}^n$ for all $\rho\geq0$, we conclude that $u\equiv 1$. Then $(\Sigma,\gamma)$ is a smooth embedded compact hypersurface in the Euclidean space with constant mean curvature and so it is a round sphere. The converse is obvious since for a Euclidean round sphere with radius $r>0$, its first Dirac eigenvalue equals to $(n-1)/(2r)$ and its mean curvature equals to $(n-1)/r$. 

If now we assume that the hypersurface $\Sigma$ is star-shaped with positive scalar and mean curvatures, we apply the method developed by Eichmair, Miao and Wang \cite{EichmairMiaoWang}. Under these assumptions, there exists a smooth map $F:\Sigma\times[0,\infty)\rightarrow \mathbb{R}^n$ such that 
\begin{eqnarray*}
\frac{\partial F}{\partial \rho}=\frac{\mathcal{H}_\rho}{\Ss_\rho}\nu,
\end{eqnarray*}
$F(\Sigma,0)=\Sigma$ and $\Sigma_\rho:=F(\Sigma,\rho)$ has positive mean and scalar curvatures and is a strictly convex hypersurface for $\rho$ sufficiently large, say for $\rho\geq \rho_0$. Here $\nu$ is the outer unit normal to the hypersurface $\Sigma_\rho$, $\mathcal{H}_\rho$  is its mean curvature and $\Ss_\rho$ its scalar curvature. Moreover, the pull-back of the Euclidean metric by $F$ on $\Sigma\times[0,\rho_0]$ has the form $g_\eta$ given by (\ref{QS-metric}) with $\eta$ a smooth positive function. From \cite[Proposition 2]{EichmairMiaoWang}, there exists a smooth positive function $v$ on $\Sigma\times [0,\rho_0]$ such that the scalar curvature of the metric $g_v$ 
is zero and with initial data
\begin{eqnarray*}
v(\cdot,0)=\frac{\eta(\cdot,0)}{2\lambda_1(\Sigma,\gamma)}\mathcal{H}_0.
\end{eqnarray*} 
This last condition implies that the mean curvature of $\Sigma\simeq\Sigma\times\{0\}$ equals to $h_0=2\lambda_1(\Sigma,\gamma)$ for the metric $g_v$. Then, as proved in \cite[Proposition 3]{EichmairMiaoWang}, the function
\begin{eqnarray*}
\rho\in[0,\rho_0]\mapsto\int_{\Sigma_\rho}\big(\mathcal{H}_\rho-h_\rho\big)d\Sigma_\rho
\end{eqnarray*}
where $h_\rho$ represents the mean curvature of $\Sigma_\rho$ with respect to the metric $g_v$, is monotone non-increasing in $\rho$ and so it follows that 
\begin{eqnarray}\label{Monotone1}
\int_{\Sigma}\big(\mathcal{H}_0-h_0\big)d\Sigma=\int_{\Sigma}\mathcal{H}_0d\Sigma-2\lambda_1(\Sigma,\gamma)|\Sigma| \geq \int_{\Sigma}(\mathcal{H}_{\rho_0}-h_{\rho_0})d\Sigma.
\end{eqnarray}
Now we can apply the original construction of Shi and Tam to the strictly convex hypersurface $\Sigma_{\rho_0}$ by solving the quasi-spherical equation in such a way that the resulting manifold $(\Sigma\times[\rho_0,\infty),g_u)$ is an asymptotically flat manifold whose inner boundary is isometric to $(\Sigma,\gamma_{\rho_0})$ with mean curvature equals to $h_{\rho_0}$. As before, we get from (\ref{GeneralInequality}) that
\begin{eqnarray}\label{Monotone2}
\int_{\Sigma}(\mathcal{H}_{\rho_0}-h_{\rho_0})d\Sigma\geq c(n)m_{ADM}(\mathcal{E},g_u).
\end{eqnarray}
Gluing $(\Sigma\times[0,\rho_0],g_v)$ and $(\Sigma\times[\rho_0,\infty),g_u)$ along their common boundary $(\Sigma,\gamma_{\rho_0})$ yields an asymptotically flat manifold $(\mathcal{E},\widetilde{g})$ with Lipschitz metric along $\Sigma_{\rho_0}$, with ADM mass equals to $m_{ADM}(\mathcal{E},g_u)$ and with a compact inner boundary whose mean curvature equals to $2\lambda_1(\Sigma,\gamma)$. Combining (\ref{Monotone1}) and (\ref{Monotone2}) with the fact that Theorem \ref{HerzlichPMT-n} holds in this context (see Remark \ref{SpinBoundaryCorner}) yields the desired inequality. If equality occurs, we deduce as before that the manifold $(\mathcal{E},\widetilde{g})$ is flat away from $\Sigma_{\rho_0}$ and that $u\equiv 1$ since $\Sigma_\rho$ is convex for $\rho\geq \rho_0$. On the other hand, we also deduce from the Gauss formula that $u\equiv \eta$ since $\Sigma_\rho$ has positive scalar curvature for $\rho\in[0,\rho_0]$. This allows to conclude that $(\mathcal{E},\widetilde{g})$ is isometric to the exterior of $\Sigma$ in the Euclidean space. The end of the proof then proceeds as in the previous case. 
\hfill$\square$

\begin{remark}
The Cauchy-Schwarz inequality gives
\begin{eqnarray*}
\int_\Sigma \mathcal{H}_0d\Sigma\leq |\Sigma|^{1/2}\Big(\int_\Sigma \mathcal{H}_0^2d\Sigma\Big)^{1/2}
\end{eqnarray*}
so that our inequality (\ref{NewUpperBound1}) implies the well-known upper bound  
\begin{eqnarray*}\label{BarIneq}
\lambda_1(\Sigma,\gamma)^2\leq\frac{1}{4}
\frac{\int_\Sigma \mathcal{H}_0^2d\Sigma}{|\Sigma|}
\end{eqnarray*}
due to B\"ar \cite{Bar1}. Note however that this last inequality holds in a much more broader context since it only assumes the existence of an isometric immersion of $(\Sigma,\gamma)$ in a Riemannian manifold carrying a parallel spinor field. 
\end{remark}

\begin{remark}\label{MinkowskiInequality-EuclideanSpace}
When $\Sigma$ is homeomorphic to a $2$-sphere, the first eigenvalue of the Dirac operator on $(\Sigma,\gamma)$ satisfies the B\"ar-Hijazi inequality \cite{Hijazi2,Hijazi1,Bar3}  
\begin{eqnarray}\label{BarHijazi}
\lambda_1(\Sigma,\gamma)\geq 2\sqrt{\frac{\pi}{|\Sigma|}}
\end{eqnarray} 
with equality if, and only if, $(\Sigma,\gamma)$ is isometric to a round sphere. If now we assume that $(\Sigma,\gamma)$ satisfies the assumptions of Theorem \ref{NewUpperBound}, one can apply (\ref{NewUpperBound1}) for $n=3$ which, with (\ref{BarHijazi}), yields the well-known Minkowski inequality 
\begin{eqnarray*}
\int_\Sigma \mathcal{H}_0d\Sigma\geq 4\sqrt{\pi\,|\Sigma|}.
\end{eqnarray*}
The equality holds only for the round spheres.
\end{remark}

\begin{remark}\label{SpinBoundaryCorner}
We briefly explain here how to prove Theorem \ref{HerzlichPMT-n} for asymptotically flat manifolds with corners along a hypersurface. An asymptotically flat manifold $(M,g)$ with boundary $\Sigma$ is said to have a corner along a hypersurface $N$ if it can be written as the disjoint union of two subsets $M_+$ and $M_-$ where $(M_+,g_{+})$ is a smooth $n$-dimensional asymptotically flat manifold with boundary $N$ and $(M_-,g_-)$ is a smooth $n$-dimensional compact Riemannian manifold with boundary $\Sigma\cup N$ and where $g_\pm=g_{|M_\pm}$. The ADM mass of $(M,g)$ corresponds to the ADM mass of $(M_+,g_+)$. Let also $H_+$ and $H_-$ the mean curvatures of $N$ in $(M_+,g_+)$ and $(M_-,g_-)$ with respect to the unit normal pointing toward infinity. The positive mass theorem for such manifolds asserts that if the scalar curvature of $(M,g)$ is nonnegative and if $H_-\geq H_+$ then its ADM mass is nonnegative. The proof of the positive mass theorem \ref{HerzlichPMT-n} relies essentially on two important points: the existence of a harmonic spinor field satisfying the Atiyah-Patodi-Singer condition which is asymptotic to a constant spinor, usually called a Witten spinor, and an integral version of the Schr\"odinger-Lichnerowicz formula. For the first point, it is enough to observe that since the metric $g$ is Lipschitz, the work of Bartnik and Chru\'sciel \cite{BartnikChrusciel} still applies to construct an adequate spinor $\Phi$. On the other hand, recall that the Schr\"odinger-Lichnerowicz formula asserts that
\begin{eqnarray*}\label{SL-Formula}
D^2\psi=\nabla^*\nabla\psi+\frac{R}{4}\psi
\end{eqnarray*}
for all spinor field $\psi$ on $M$ and where $D$, $\nabla$ and $\nabla^\ast$ are respectively the Dirac operator, the spin Levi-Civita connection and its $L^2$-formal adjoint on $(M,g)$. Although this formula is no longer valid on $(M,g)$ since the scalar curvature is not defined on $N$, it remains true on $(M_+,g_+)$ and on $(M_-,g_-)$ separately. We can then integrate this formula on $M_-$ first and then on a domain $M_R$ in $M_+$ whose boundary is the union of $N$ and a coordinate sphere at infinity with radius $R$. Adding them together and inserting the aforementioned Witten spinor $\Phi$ lead to 
\begin{eqnarray*}
\frac{1}{2}(n-1)\omega_{n-1}m_{ADM}(M,g) & = & \int_{M}\Big(|\nabla\Phi|^2+\frac{R}{4}|\Phi|^2\Big)\,dM+\frac{1}{2}\int_N(H_--H_+)|\Phi|^2\,dN
\\ & & -\int_{\Sigma}\<\D_\gamma\Phi+\frac{H}{2}\Phi,\Phi\>d\Sigma
\end{eqnarray*}
as $R$ goes to infinity. The nonnegativity of the mass then follows from the assumptions $R\geq 0$, $H_-\geq H_+$ and (\ref{DiracMeanCur-n}). 
\end{remark}

 
\section{Proof of Theorem \ref{AsymptoticDiracSphere}}\label{LargeSphereLimit}


If $r$ is large enough, the Gauss curvature of $S_r$ endowed with the metric $g_r$ is positive. So, by the solution of the Weyl's embedding problem, $(S_r,g_r)$ can be isometrically embedded in $\mathbb{R}^3$ with positive mean curvature $\mathcal{H}_0$. This embedding is unique up to an isometry of $\mathbb{R}^3$. It follows from \cite[Lemma 2.4]{FanShiTam} that 
\begin{eqnarray}\label{TotalMC-Euclidean}
\int_{S_r}\mathcal{H}_0dS_r=4\pi r+\frac{|S_r|}{r}+O(r^{1-2\tau})
\end{eqnarray}
which, when combined with Corollary \ref{2SphereDirac}, yields
\begin{eqnarray*}
\lambda_1(S_r,g_r)|S_r|\leq 2\pi r+\frac{|S_r|}{2r}+O(r^{1-2\tau}).
\end{eqnarray*}
From \cite[Lemma 2.1]{FanShiTam}, the area of $S_r$ has the following expansion 
\begin{eqnarray}\label{AreaExpansion}
|S_r|=4\pi r^2+\beta(r)+O(r^{1-2\tau})
\end{eqnarray}
near infinity where 
\begin{eqnarray*}
\beta(r)=\frac{1}{2}\int_{S_r}g_r^{ij}\sigma_{ij}dS_r
\end{eqnarray*}
satisfies $\beta(r)=O(r^{2-\tau})$ so that
\begin{eqnarray}\label{UpperBoundExp}
\lambda_1(S_r,g_r)|S_r|\leq 4\pi r+\frac{\beta(r)}{2r}+O(r^{1-2\tau}).
\end{eqnarray}
Here $g_r^{ij}$ is $(g_r)_{ij}$ raised with respect to $g^{ij}$. On the other hand, the B\"ar-Hijazi inequality (\ref{BarHijazi}) with (\ref{AreaExpansion}) leads to 
\begin{eqnarray}\label{LowerBoundExp}
\lambda_1(S_r,g_r)|S_r|\geq 2\sqrt{\pi}|S_r|^{1/2}=4\pi r+\frac{\beta(r)}{2r}+O(r^{1-2\tau}).
\end{eqnarray}
From (\ref{UpperBoundExp}) and (\ref{LowerBoundExp}) we deduce that
\begin{eqnarray}\label{DiracExpansion}
\lambda_1(S_r,g_r)|S_r|=4\pi r+\frac{\beta(r)}{2r}+o(1)
\end{eqnarray}
since $\tau>\frac{1}{2}$. Now recall from \cite[Lemma 2.2]{FanShiTam} that
\begin{eqnarray*}
\int_{S_r}H_rdS_r=\frac{|S_r|}{r}+4\pi r-8\pi m_{ADM}(M,g)+o(1)
\end{eqnarray*}
which, with (\ref{AreaExpansion}), gives 
\begin{eqnarray}\label{MeanCurvatureExpansion}
\frac{1}{2}\int_{S_r}H_rdS_r= 4\pi r+\frac{\beta(r)}{2r}-4\pi m_{ADM}(M,g)+o(1).
\end{eqnarray}
The asymptotic expansion (\ref{IntegralInequality1}) follows directly from (\ref{DiracExpansion}) and (\ref{MeanCurvatureExpansion}). 
\hfill$\square$

\begin{remark}
In \cite{HijaziMontielZhang1}, Hijazi, Montiel and Zhang proved an inequality relating the first nonnegative eigenvalue of the Dirac operator of hypersurfaces bounding compact spin Riemannian manifolds with nonnegative scalar curvature. A direct consequence of this result ensures that if $(\Omega,g)$ is a $3$-dimensional manifold whose boundary $\Sigma$ is the union of a minimal $2$-sphere $\Sigma_H$ and a surface $\Sigma_O$ with positive mean curvature $H_O$, the first nonnegative eigenvalue $\lambda_1(\Sigma_O,\gamma)$ of the Dirac operator $(\Sigma_O,\gamma)$ with $\gamma:=g_{|\Sigma_O}$ satisfies
\begin{eqnarray}\label{HMZ-Extrinsic}
\lambda_1(\Sigma_O,\gamma)>\frac{1}{2}\min_{\Sigma} H_O.
\end{eqnarray}
The proof of the inequality (\ref{HMZ-Extrinsic}), in a broad sense, is given in \cite[Theorem 6]{HijaziMontielZhang1}. The fact that the equality cannot hold for a domain $\Omega$ as above can be seen as follows. In fact, it is easy to see that if equality holds in (\ref{HMZ-Extrinsic}), the domain $\Omega$ carries a parallel spinor with respect to the metric $g$ whose restriction to $\Sigma_H$ gives rise to a harmonic spinor on this boundary component. However, this is impossible since we assumed that $\Sigma_H$ is a $2$-sphere on which such a spinor field cannot exist by the B\"ar-Hijazi inequality (\ref{BarHijazi}). This inequality applies for example when $(M,g)$ is a $3$-dimensional asymptotically flat manifold with nonnegative scalar curvature and with compact minimal inner boundary $\partial M$ and $\Omega_r$ is the compact domain whose boundary is the disjoint union of $\partial M$ and a coordinate sphere $S_r$ with $r>0$. The asymptotic expansion (\ref{IntegralInequality1}) shows explicitly that the non-sharpness of (\ref{HMZ-Extrinsic}) in this situation is directly related to the positivity of the ADM mass of $(M,g)$. Indeed, from the work of Bray \cite[Theorem 9]{Bray}, the presence of a minimal compact boundary on $(M,g)$ ensures that its ADM mass is positive and so it follows from Theorem \ref{AsymptoticDiracSphere} that 
\begin{eqnarray}\label{IntegralLowerBound}
\lambda_1(S_r,g_r)>\frac{1}{2|S_r|}\int_{S_r}H_rdS_r\geq \frac{1}{2}\min_{S_r}H_r
\end{eqnarray}
for $r$ sufficiently large. Note that (\ref{IntegralLowerBound}) improves (\ref{HMZ-Extrinsic}) since it replaces a pointwise bound by an integral one.
\end{remark}

\begin{remark}
Assume that $(M,g)$ is a $3$-dimensional asymptotically Schwarzschild manifold of mass $m$ that is there exists a bounded set $K$ such that $M\setminus K$ is diffeomorphic to the complement of a closed ball in $\mathbb{R}^3$ and such that, in the coordinate charts induced by this diffeomorphism, the metric satisfies
\begin{eqnarray*}
g_{ij}=\Big(1+\frac{2m}{|x|}\Big)\delta_{ij}+\sigma_{ij},\quad \sum_{l=0}^4|x|^{l}|\partial^l\sigma_{ij}|=O(|x|^{-2})
\end{eqnarray*}
where $m$ is a real number. Note that the parameter $m$ corresponds exactly to the ADM mass of $(M,g)$. In this situation, one can say more since, using the same method and the estimates in \cite[Section 5]{ShiTam1}, we get an asymptotic expansion for the first eigenvalue, namely
\begin{eqnarray*}\label{ExpansionDiracEigenvalue}
\lambda_1(S_r,g_r)=\frac{1}{r}-\frac{m}{r^2}+O\Big(\frac{1}{r^3}\Big).
\end{eqnarray*} 
A direct consequence of this expansion and the positive mass theorem is the following comparison result for the first eigenvalue of the Dirac operator. If the scalar curvature of $(M,g)$ is nonnegative, then 
\begin{eqnarray*}
\lim_{r\rightarrow\infty} r^2\big(\lambda_1(S_r,\delta_r)-\lambda_1(S_r,g_r)\big)\geq 0
\end{eqnarray*}
where $\lambda(S_r,\delta_r)=1/r$ is the first eigenvalue of the Dirac operator $\D_{\delta_r}$ on the Euclidean sphere of radius $r$. Moreover, equality occurs if, and only if, $(M,g)$ is isometric to the Euclidean space. In particular, if the mass $m$ of $(M,g)$ is positive, then $\lambda_1(S_r,g_r)<\lambda_1(S_r,\delta_r)$ for sufficiently large $r$. This is the case for example when $(M,g)$ has nonnegative scalar curvature and a compact inner boundary with nonpositive mean curvature. 
\end{remark}


\section{Proof of Theorem \ref{SmallSphereExpansion}}\label{SmallSphereLimit}


In this section, we give a lower and an upper estimates for the first eigenvalue of the Dirac operator which imply the expansion of Theorem \ref{SmallSphereExpansion}. Let $(M,g)$ be a $3$-dimensional Riemannian manifold and let $p\in M$ be an interior point. Consider $(x_1,x_2,x_3)$ the normal coordinates near $p$ and let $r$ be the geodesic distance form $p$. Lemma $3.2$ in \cite{FanShiTam} ensures that 
\begin{eqnarray}\label{AreaNormal}
|S_r|=4\pi r^2-\frac{2\pi}{9}R(p)r^4+\frac{\pi}{675}\big(4R(p)^2-2|Ric(p)|^2-9\Delta R(p)\big)r^6+O(r^7)
\end{eqnarray}
and so it follows from the B\"ar-Hijazi inequality (\ref{BarHijazi}) that
\begin{eqnarray}\label{SSLower}
\lambda_1(S_r,\gamma_r)\geq\frac{1}{r}+\frac{R(p)}{36}r+\frac{L(p)}{5400}r^3+O(r^4)
\end{eqnarray}
where we let
\begin{eqnarray*}
L(p)=\frac{9}{4}R^2(p)+2|Ric|^2(p)+9\Delta R(p).
\end{eqnarray*}
Here $|Ric|$ denotes the norm of the Ricci curvature and $\Delta$ is the Laplacian of $(M,g)$. On the other hand, for $r$ small enough, the sphere $(S_r,\gamma_r)$ has positive Gauss curvature in such a way that it can be isometrically embedded in $\mathbb{R}^3$ with positive mean curvature $\mathcal{H}_0$. Then it is proved in \cite[p. 66]{FanShiTam} that
\begin{eqnarray*}
\int_{S_r}\mathcal{H}_0dS_r=8\pi r-\frac{2\pi}{9}R(p)r^3-\frac{\pi}{2700}\big(99R^2(p)-312|Ric|^2(p)+36\Delta R(p)\big)r^5+O(r^6).
\end{eqnarray*}
Combining this formula with (\ref{AreaNormal}) and (\ref{NewUpperBound1}) yields  
\begin{eqnarray}\label{SSUpper}
\lambda_1(S_r,\gamma_r)\leq\frac{1}{r}+\frac{R(p)}{36}r+\frac{1}{5400}\big(L(p)+80\big|E|^2(p)\big)r^3+O(r^4)
\end{eqnarray}
where $E:=Ric-(R/3)g$ is the traceless part of the Ricci tensor of $(M,g)$. It is now obvious to deduce Theorem \ref{SmallSphereExpansion} and Corollary \ref{SmallSphereDiracV} from (\ref{SSLower}) and (\ref{SSUpper}). 

\begin{remark}
It follows from Corollary \ref{SmallSphereDiracV} that, for $r$ small enough, the integral bound 
\begin{eqnarray*}\label{IntegralLowerBoundSmallSpheres}
\lambda_1(S_r,\gamma_r)>\frac{1}{2|S_r|}\int_{S_r}H_rdS_r
\end{eqnarray*}
holds on any geodesic spheres $S_r$ centered at an interior point $p$ of any $3$-dimensional Riemannian manifolds as soon as $R(p)>0$. Once again this inequality improves (\ref{HMZ-Extrinsic}) in this situation.
\end{remark}


\section{A few words on the hyperbolic setting}\label{HyperbolicPart}


In this last section, we give the proofs of Theorem \ref{NewUpperBoundHyperbolique} and Corollary \ref{MinkowskiHyperbolic}. Let us briefly recall the setting and the main results which are needed for this purpose.

When $(\Sigma,\gamma)$ is a compact hypersurface isometrically embedded in an $n$-dimensional spin Riemannian manifold $(M,g)$ with scalar curvature bounded from below by $-n(n-1)\kappa^2$, $\kappa>0$, it is natural to consider the operators given by
\begin{eqnarray*}
\D_\gamma^\pm:=\D_\gamma\pm \frac{n-1}{2}\kappa \sqrt{-1}c(\nu)
\end{eqnarray*}
where $c(\nu)$ represents the Clifford multiplication, with respect to $g$, by the inner unit normal to $\Sigma$ denoted by $\nu$. These are first order elliptic and self-adjoint differential operators which acts on the spinor bundle over $(M,g)$ restricted to $\Sigma$. They appear as natural counterparts of the extrinsic Dirac operator in the integral version of the hyperbolic Schr\"odinger-Lichnerowicz formula (see for example \cite{ChruscielHerzlich,WangYau0,Kwong1,HijaziMontielRoldan,HijaziMontielRaulot3}). Since the Clifford multiplication by $\nu$ sends an eigenspinor for $\D^+_\gamma$ associated to $\lambda$ to an eigenspinor for $\D^-_\gamma$ associated to $-\lambda$, the spectra of these operators, denoted by $Spec(\D^\pm_\gamma)$, are such that 
\begin{eqnarray}\label{SpectreRelation}
Spec(\D_\gamma^+)=-Spec(\D_\gamma^-)\subset\mathbb{R}^*.
\end{eqnarray}
Let $\lambda^\pm_1(\Sigma,\gamma)$ denote their first nonnegative eigenvalues. Note that since 
\begin{eqnarray*}
\big(\D_\gamma^\pm)^2\psi=\D_\gamma^2\psi+\frac{(n-1)^2}{4}\kappa^2\psi
\end{eqnarray*}
for all spinor fields $\psi$ on $M$ restricted to $\Sigma$, it is direct to deduce that 
\begin{eqnarray}\label{SpectreDpm}
\lambda^\pm_1(\Sigma,\gamma)^2\geq \lambda_1(\Sigma,\gamma)^2+\frac{(n-1)^2}{4}\kappa^2.
\end{eqnarray}
Let us now prove Theorem \ref{NewUpperBoundHyperbolique}.

\vspace{0.2cm}

{\it Proof of Theorem \ref{NewUpperBoundHyperbolique}:} We proceed exactly as in the proof of Theorem \ref{NewUpperBound}. From \cite{WangYau0,Kwong1}, it follows that under our assumptions, there exists on $\mathbb{H}^n(-\kappa^2)\setminus\Omega\simeq\Sigma\times[0,+\infty[$ an unique function $w$ with initial value
\begin{eqnarray*}
w(\cdot,0)=\frac{1}{2\lambda^\pm_1(\Sigma,\gamma)}\mathcal{H}_0
\end{eqnarray*}
in such a way that the quasi-spherical metric $g_w:=w^2d\rho^2+\gamma_\rho$ defines an asymptotically hyperbolic metric with constant scalar curvature $R=-n(n-1)k^2$ and with $(\Sigma,\gamma)$ as a compact inner boundary with mean curvature equals to $H=2\lambda_1^\pm(\Sigma,\gamma)$. Here $\Omega$ denotes the compact domain of $\mathbb{H}^n(-\kappa^2)$ bounded by $\Sigma$ and $\rho$ is the distance from $\Sigma$. Then it follows from the positive mass theorem recalled in Remark \ref{PMT-HyperbolicBoundary} that 
\begin{eqnarray}\label{MassHyperQS}
\lim_{\rho\rightarrow+\infty}\int_{\Sigma_\rho}\mathcal{H}_{\rho}(1-w^{-1}) X\cdot\zeta\,d\Sigma_\rho \leq 0
\end{eqnarray}
for any future-directed null vector $\zeta\in\mathbb{R}^{n,1}$. The vector $X=(x_1,\cdots,x_n,t)$ is the position vector in $\mathbb{R}^{n,1}$, the inner product is given by the Lorentz metric and $\mathcal{H}_\rho$ denotes the mean curvature of $\Sigma_\rho$ in $\mathbb{H}^n(-\kappa^2)$. This implies in particular that 
\begin{eqnarray}\label{MassHyperQS-1}
\lim_{\rho\rightarrow+\infty}\int_{\Sigma_\rho}\mathcal{H}_{\rho}(1-w^{-1}) \cosh(\kappa r)\,d\Sigma_\rho \geq 0.
\end{eqnarray}
On the other hand, it is proved in \cite{ShiTam2,Kwong1} that there exists $\alpha>1$ such that for any future-directed null vector $\zeta\in\mathbb{R}^{n,1}$ the function
\begin{eqnarray*}\label{MonotonicityFormula}
\rho\in[0,+\infty[\mapsto \int_{\Sigma_\rho}\mathcal{H}_{\rho}(1-w^{-1})X_\alpha\cdot\zeta\,d\Sigma_\rho
\end{eqnarray*} 
is nonincreasing in $\rho$ where $X_\alpha=(x_1,\cdots,x_n,\alpha t)$. Combining this fact with (\ref{MassHyperQS}) and (\ref{MassHyperQS-1}) yields
\begin{eqnarray*}
\int_{\Sigma}\big(\mathcal{H}_{0}-2\lambda_1^\pm(\Sigma,\gamma)\big)X_\alpha\cdot\zeta\,d\Sigma\leq\lim_{\rho\rightarrow+\infty}\int_{\Sigma_\rho}\mathcal{H}_{\rho}(1-w^{-1}) X_\alpha\cdot\zeta\,d\Sigma_\rho \leq 0
\end{eqnarray*}
and this implies that the vector 
\begin{eqnarray*}
\int_{\Sigma}\big(\mathcal{H}_{0}-2\lambda_1^\pm(\Sigma,\gamma)\big)X_\alpha\,d\Sigma
\end{eqnarray*}
is future-directed causal. In particular, its time component is nonnegative and so the inequality (\ref{NewUpperBoundHyp1}) follows straightforwardly. The equality case is a direct consequence of the characterization of the equality case in the positive mass theorem given in Remark \ref{PMT-HyperbolicBoundary}. 
\hfill$\square$

\vspace{0.2cm}

{\it Proof of Corollary \ref{MinkowskiHyperbolic}:} The B\"ar-Hijazi inequality (\ref{BarHijazi}) and the estimate (\ref{SpectreDpm}) ensure that 
\begin{eqnarray*}
\lambda^\pm_1(\Sigma,\gamma)^2\geq\frac{4\pi}{|\Sigma|}+\kappa^2. 
\end{eqnarray*}
Moreover, equality holds only for two dimensional round spheres. Then the Minkowski-type inequality (\ref{MinkowskiHyperbolic}) as well as its equality case follow directly from Theorem \ref{NewUpperBoundHyperbolique}.
\hfill$\square$

\begin{remark}
By sending $\kappa\rightarrow 0$ in the inequalities (\ref{NewUpperBoundHyp1}) and (\ref{MinkowskiHyperbolic1}), we respectively recover the estimate (\ref{NewUpperBound}) and the classical Minkowski inequality for convex body in the $3$-dimensional Euclidean space (see Remark \ref{MinkowskiInequality-EuclideanSpace}). 
\end{remark}

\begin{remark}
It is obvious to deduce from the inequality (\ref{NewUpperBoundHyp1}) that 
\begin{eqnarray*}
\lambda_1^\pm(\Sigma,\gamma)\leq\frac{1}{2}\sup_\Sigma\mathcal{H}_0
\end{eqnarray*}
which imply, with the help of (\ref{SpectreDpm}), a well-known inequality due to Ginoux \cite{Ginoux2} and which asserts that 
\begin{eqnarray*}
\lambda_1(\Sigma,\gamma)^2\leq\frac{1}{4}\big((\sup_\Sigma\mathcal{H}_0)^2-(n-1)^2\kappa^2\big).
\end{eqnarray*}
Note however that this last inequality holds in a much more broader context since it only assumes the existence of an isometric immersion of $(\Sigma,\gamma)$ in a Riemannian manifold carrying an imaginary Killing spinor field. 
\end{remark}

\begin{remark}\label{PMT-HyperbolicBoundary}
In the proof of Theorem \ref{NewUpperBoundHyperbolique}, we make appeal to a positive mass theorem for asymptotically hyperbolic manifolds with compact inner boundary which we now explain. This result is inspired by a work with Hijazi and Montiel \cite{HijaziMontielRaulot4} and which can be seen as a hyperbolic version of the results of Herzlich \cite{Herzlich1,Herzlich2}. Assume, as in the proof of Theorem \ref{NewUpperBoundHyperbolique}, that $M:=\mathbb{H}^n(-\kappa^2)\setminus\Omega$ is endowed with respect to the metric $g_w$. Thus $(M,g_w)$ is an $n$-dimensional asymptotically hyperbolic manifold, in the sense of \cite{AnderssonDahl1}, with constant scalar curvature and with compact inner boundary $(\Sigma,\gamma)$ whose mean curvature is constant equals to $2\lambda^+_1(\Sigma,\gamma)$. It is obvious to see that $M$ is endowed with a spin structure and so we can consider $(SM,\nabla,c,\<\,,\,\>)$ the associated Dirac bundle for the metric $g_w$ where $SM$ denotes the bundle of complex spinors, $\nabla$ is the corresponding spin Levi-Civita connection, $c$ the Clifford multiplication and $\<\,,\,\>$ the Hermitian inner product.
On the other hand, on $(M,\mathfrak{h}_\kappa)$, $\mathfrak{h}_\kappa$ being the hyperbolic metric with constant sectional curvature $-\kappa^2$, there exists a set of maximal dimension of imaginary Killing spinors which is parametrized by $a\in\mathbb{C}^{2^m}$ with $m=[n/2]$. This implies that for every such a $a\in\mathbb{C}^{2^m}$ corresponds an imaginary Killing spinor field $\phi_a$ which can be considered as a section of $SM$ via the identification between the spinor bundles over $(M,g_w)$ and $(M,\mathfrak{h}_\kappa)$. Then it can be shown that there exists an unique $\psi_a\in W^{1,2}$ such that the spinor field $\Psi_a:=\psi_a+\eta\Phi_a\in\Gamma(SM)$ satisfies the following boundary value problem:
$$\left\lbrace
\begin{array}{ll}
D^+\Psi_a=0 \quad& {\rm on } \,\,M,\nonumber\\
\Paps^+\Psi_{a\,|\Sigma}=0 \quad& {\rm along }\,\,\Sigma\nonumber.
\end{array}
\right.$$
Here $\eta$ is a cut-off function that vanishes on a compact set of $M$ and is equal to $1$ for $r$ large enough. Moreover, $D^\pm:=c\circ\nabla^\pm$ is a zero order modification of the classical Dirac operator $D$ with $\nabla^\pm$ the modified connection defined by 
\begin{eqnarray*}
\nabla^\pm_Z=\nabla_Z\pm\frac{\sqrt{-1}}{2}\kappa c(Z)
\end{eqnarray*}
for all $Z\in\Gamma(TM)$. On the other hand, the maps $\Paps^\pm$ represent the $L^2$-orthogonal projections onto the subspace spanned by the eigenspinors corresponding to the positive eigenvalues of $\D_\gamma^\pm$ which define global elliptic boundary conditions for the operators $D^\pm$. Then integrating by parts on the compact domain delimited by $\Sigma_\rho$ the identity
\begin{eqnarray*}
(\nabla^+)^*\nabla^+\Psi_a=0
\end{eqnarray*}
which is deduced from the Schr\"odinger-Lichnerowicz formula and sending $\rho\rightarrow+\infty$ finally leads to 
\begin{eqnarray}\label{MassBoundary}
\lim_{\rho\rightarrow+\infty}\int_{\Sigma_\rho}\mathcal{H}_r(1-w^{-1})|\Phi_a|^2\,d\Sigma_{\rho} =\int_M|\nabla^+\Psi_a|^2\,dM-\int_\Sigma\<\D^+_\gamma\Psi_a+\frac{1}{2}H\Psi_a,\Psi_a\>\,d\Sigma
\end{eqnarray}
which is nonnegative since $\Paps^+\Psi_{a|\Sigma}=0$ and $H=2\lambda^+_1(\Sigma,\gamma)$. The fact that the left-hand side of the previous equality is finite is proved in \cite{WangYau0,Kwong1}. Now the inequality (\ref{MassHyperQS}) follows directly since for every null vector $\zeta\in\mathbb{R}^{n,1}$, there exists $a\in\mathbb{C}^{2^m}$ where 
\begin{eqnarray*}
\zeta_a=\sum_{j=1}^n\<\sqrt{-1}c_0(e_j)a,a\>e_j-\<\sqrt{-1}c_0(e_0)a,a\>e_0
\end{eqnarray*}
such that
\begin{eqnarray*}
-2\kappa X\cdot\zeta_a=|\Phi_a|^2.
\end{eqnarray*}
Here $e_0=\frac{\partial}{\partial t}$, $e_j=\frac{\partial}{\partial x_j}$ in $\mathbb{R}^{n,1}$ for all $j=1,\cdots,n$ and $c_0$ denotes the Clifford multiplication with respect to the Lorentz inner product in $\mathbb{R}^{n,1}$. Finally note that if the left-hand side of (\ref{MassBoundary}) vanishes, the spinor $\Psi_a$ has to be an imaginary Killing spinor and so in particular $(M,g_w)$ is an Einstein manifold with scalar curvature equals to $-n(n-1)\kappa^2$. Then, from the Gauss equation for the embedding of $\Sigma_\rho$ in $(M,g_w)$ we deduce that
\begin{eqnarray*}
(1-w^{-2})(\Ss_\rho+(n-1)(n-2)\kappa^2)=0
\end{eqnarray*}
which implies that $w\equiv 1$ since we assumed that the sectional curvature $K$ of $(\Sigma,\gamma)$ satisfies $K>-\kappa^2$.  This implies that $\Sigma$, as inner boundary of $M$, has the same second fundamental form as boundary of the hyperbolic domain $\Omega$. Then $(\Sigma,\gamma)$ is a smooth embedded compact hypersurface in the hyperbolic space with constant mean curvature and so it is a round sphere. 
\end{remark}


\bibliographystyle{alpha}   
\bibliography{BiblioHabilitation}     


\end{document}